\documentclass{article}
\usepackage{amsmath}
\usepackage{amssymb}
\usepackage[dvips]{epsfig}
\usepackage[arrow,matrix]{xy}
\newcommand{\Om}{\Omega}
\newcommand{\Oms}{\Omega^{\mathrm{Spin}}}
\newcommand{\fN}{\mathfrak N}
\newcommand{\bR}{\mathbb R}

\newcommand{\bP}{\mathbb P}
\newcommand{\cC}{\mathcal C}

\newcommand{\bZ}{\mathbb Z}
\newcommand{\bF}{\mathbb F}
\newcommand{\Tor}{\operatorname{Tor}}

\newcommand{\co}{\colon\,}

\newcommand{\lp}{\textup{(}}
\newcommand{\rp}{\textup{)}}
\newcommand{\scurv}{scalar curvature}
\newcommand{\psc}{positive scalar curvature}
\newcommand{\Yi}{Yamabe invariant}
\newcommand{\nYi}{nonnegative Yamabe invariant}
\newcommand{\dvol}{d\mbox{\textup{vol}}}
\newcommand{\diam}{\operatorname{diam}}
\newcommand{\vol}{\operatorname{vol}}
\newcommand{\tor}{^{\text{\textup{toral}}}}
\newcommand{\ator}{^{\text{\textup{atoral}}}}
\newtheorem{thm}{Theorem}[section]
\newtheorem{cor}[thm]{Corollary}
\newtheorem{lemma}[thm]{Lemma}
\newtheorem{prop}[thm]{Proposition}
\newtheorem{remk}[thm]{Remark}
\newtheorem{defin}[thm]{Definition}
\newtheorem{examples}[thm]{Examples}
\newtheorem{prob}[thm]{Problem}
\newenvironment{demo}{\par\noindent{\em
Proof}. }{$\square$\par\indent}
\newenvironment{sketch}{\par\noindent{\em
Sketch of Proof}. }{$\square$\par\indent}
\newenvironment{labeleddemo}[1]{\par\noindent{\em
#1}. }{$\square$\par\indent}
\newenvironment{defn}{\begin{defin}\em}{\end{defin}}
\newenvironment{exs}{\begin{examples}\em}{\end{examples}}

\numberwithin{equation}{section}
\begin{document}
\title{The Yamabe invariant for\\ non-simply connected manifolds}
\author{Boris Botvinnik and Jonathan Rosenberg}
\maketitle
\maketitle
\markright{The non-simply connected Yamabe invariant}
\begin{abstract}
The Yamabe invariant is an invariant of a closed smooth manifold
defined using conformal geometry and the
scalar curvature. Recently, Petean 
showed that the Yamabe invariant is non-negative for all closed
simply connected manifolds of dimension $\ge 5$. We extend this
to show that Yamabe invariant is non-negative for all closed
manifolds of dimension $\ge 5$ with fundamental group
of odd order having all Sylow subgroups abelian.
The main new geometric input is a way of studying the
Yamabe invariant on Toda brackets.
A similar method of proof shows that all closed
manifolds of dimension $\ge 5$ with fundamental group
of odd order having all Sylow subgroups elementary
abelian, with non-spin universal
cover, admit metrics of \psc, once one restricts to the
``complement'' of manifolds whose
homology classes are ``toral.''  The exceptional 
toral homology classes only exist in dimensions not exceeding the
``rank'' of the fundamental group, so this proves important cases of the
Gromov-Lawson-Rosenberg Conjecture once the dimension is sufficiently
large. 
\end{abstract}

\section{Introduction}\label{sec:intro}

The positive solution of the Yamabe problem \cite{Sc} tells us that
if $M$ is a compact smooth manifold (without boundary), then every
conformal class $\cC$ of Riemannian metrics on $M$ contains a metric 
(known as a Yamabe metric) of constant scalar curvature with the
following special property. Its scalar curvature is the infimum of the
scalar curvature $s_g$, taken over all metrics in $\cC$ with constant
scalar curvature and total volume $1$. The value
of this scalar curvature is called the Yamabe constant
$Y(M,\cC)$ of $\cC$. Equivalently, $Y(M,\cC)$ can be defined to be the minimum
over metrics $g\in\cC$ of the Einstein-Hilbert functional
\[
I(g) =\frac{\int_M s_g \ d\vol_g}{\vol_g(M)^{\frac{n-2}{n}}}.
\]
It is therefore natural to to ask if there is a ``best'' Yamabe
metric, and if so what its scalar curvature is. 
That motivates the following definition from \cite{Kob}. The Yamabe 
invariant of $M$ is defined by 
\begin{equation}\label{eq:Ydef}
 Y(M)=\sup_\cC Y(M,\cC).
\end{equation}
This supremum is not always attained, so the answer to the question
about whether $M$ has a ``best'' metric of constant scalar curvature
might be ``no.''  The best that
is known that there are singular metrics (with singularities at a finite
number of points) which serve as the ``best'' approximation to an
Einstein metric on $M$. 

Nevertheless, $Y(M)$ is a diffeomophism invariant of $M$. It also
turns out that $Y(M)>0$ if and only if $M$ admits a metric of
positive
scalar curvature, a much-studied condition (\cite{GL1}, \cite{GL2}, 
\cite{GL3}, \cite{S1}, \cite{S2}, \cite{RS1}, \cite{RS2}, \cite{BGS}).
However, $Y(M)=0$ is possible even when $M$ admits no scalar-flat
metric.

In dimension $2$, Gauss-Bonnet quickly shows that $Y(M)=4\pi\chi(M)$.
In dimension $4$, $Y(M)$ can be positive, $0$, or negative, and a lot
is known about it from Seiberg-Witten theory (\cite{LeB} and
\cite{P1}). Similarly, there is a conjectural connection between
$Y(M)$ and ``geometrization'' when $\dim M=3$ (see for
instance \cite{And}). But even when $\dim M=3$, and especially
when $\dim M>4$, it is not yet known if
there are \textbf{any} manifolds with $Y(M)<0$. (The obvious candidates
for such manifolds are hyperbolic manifolds, but for all we know they
could have vanishing \Yi.) In fact, Petean \cite{P2}
has proved that $Y(M)\ge 0$ for any simply connected
manifold of dimension at least $5$.

In this paper we study the Yamabe invariant for manifolds with finite
fundamental groups. Our first main result is the following.
\begin{thm}\label{thm:thmA}
Let $M$ be a closed, connected, compact manifold
with finite fundamental group $\pi$, $\dim M \ge 5$. Suppose all Sylow
subgroups of $\pi$ are abelian. Assume either that $M$ is spin and
the order of $\pi$ is odd, or else that the universal cover of $M$ is
non-spin. Then $Y(M)\geq 0$.
\end{thm}
The proof of this result is somewhat involved. First of all, we use 
surgery tools (developed in the study of positive scalar curvature) to
reduce the assertion of Theorem \ref{thm:thmA} to special situations. In
particular, we show that it is enough to study the case when $\pi$ is
a finite abelian $p$-group. The central objects to understand here are
the bordism groups $\Om(B\pi)$ and $\Oms(B\pi)$, and
the proof amounts to the fact that all elements of these bordism groups
may be represented by manifolds with nonnegative Yamabe invariant. A
computation of these groups is quite hard, and the actual answer is known
only for elementary abelian groups of odd order and few other cases (see
\cite{JW}). Instead we use the K\"unneth formula to build manifolds with
nonnegative Yamabe invariant to represent generators of these bordism
groups. There are two types of ``building blocks'': tensor products
(which are realized by direct products of manifolds) and torsion
products (which geometrically are just Toda brackets). 

We recall that Toda bracket $\langle M, P, L \rangle$ is defined when
$M\times P = \partial V$ and $P\times L = \partial U$. Then the manifold
\[
W= V\times L\cup_{M\times P\times L}M\times U
\]
represents the Toda bracket $\langle M, P, L \rangle$. 

As usual, to prove new geometric results we have to employ some 
\emph{new geometric techniques}. Roughly, we show (under some
restrictions) that if $Y(M)$ and $Y(L)$ are $\ge 0$ (resp.,
$>0$), then $Y(W)\ge 0$
(resp., $>0$). We prove this by analytical means using elementary
differential geometry.

Our second main result is the following.
\begin{thm}\label{thm:thmB}
Let $M$ be a closed, connected, compact manifold
with fundamental group $\pi$ of odd order. Suppose all Sylow
subgroups of $\pi$ are elementary abelian of rank $\le r$. Assume that $M$ is
non-spin and that $\dim M \ge \max(5,r)$. Then $M$ has
a metric of {\psc}.
\end{thm}

To put these results in context, it's worth recalling what is known
about {\psc} for manifolds with finite fundamental group.
For such manifolds (of dimension $\ge 5$) whose universal cover is non-spin,
\emph{there are no known obstructions to \psc}. For spin manifolds
of dimension $\ge 5$
with finite fundamental group, the only known obstructions to
{\psc} come from the index theory of the Dirac operator
(\cite{RS1}, \cite{RS2}), and it is known that ``stably'' these
are the only obstructions \cite{RS1}. In fact in \cite{RosKO}, it
was conjectured (on the basis of extremely spotty evidence) that
the index theory of the Dirac operator provides the only
obstructions to {\psc} on manifolds of dimension $\ge 5$
with finite fundamental group. This conjecture has sometimes been
called the Gromov-Lawson-Rosenberg Conjecture. However, the ``stable''
theorem by itself does not actually answer the question of whether any 
\emph{particular} manifold with vanishing Dirac obstructions
admits a metric of {\psc}. It \emph{is} known
\cite{BGS} that for spin manifolds of  dimension $\ge 5$ with
finite fundamental group with periodic cohomology, the Dirac
obstructions are the only obstructions to \psc. A similar
theorem was proved by Schultz \cite{Sch},
and independently by Botvinnik and Gilkey \cite{BG}, for spin manifolds of  
dimension $\ge 5$ with fundamental group $\bZ/p\times\bZ/p$, $p$
an odd prime. But very little was
previously known about {\psc} for manifolds with elementary abelian
fundamental group of rank $>2$. The proof of Theorem \ref{thm:thmB}
is based on a reduction to the results of \cite{BG}, again
using Toda brackets.

The outline of the paper is as follows. Section \ref{sec:redtools}
recalls the surgery and bordism theorems necessary for attacking the
problems. Section \ref{sec:tools} contains our basic geometric results
on Toda brackets.  Section \ref{sec:Todaresults} puts together the
topological and geometrical tools to prove Theorem \ref{thm:thmA}
and related results, and Section \ref{sec:psc} proves
Theorem \ref{thm:thmB} and related results.

We would like to thank Sergey Novikov for his encouragement and support.

\section{Basic Topological Reduction Tools}\label{sec:redtools}

To warm up, we recall the following result of Petean for simply
connected manifolds:
\begin{thm}[\cite{P2}]\label{thm:Peteansimplyconn}
If $M^n$ is a connected, simply connected closed
manifold of dimension $n\ge 5$, then $Y(M)\ge0$.
\end{thm}
The proof of this fact is based on the following \emph{surgery
theorem}:
\begin{thm}[Petean, Yun \cite{PY}]\label{thm:surgery}If $M$ is a closed
manifold with connected components $M_i$, and if another closed 
connected manifold $M'$ can be obtained from
$M$ by surgeries in codimension $\ge 3$, and if $Y(M_i)\ge0$
for each $i$, then $Y(M')\ge 0$.
\end{thm}
\begin{demo}This is really three theorems in one. If $Y(M_i)>0$
for all $i$, then $M$
admits a metric of \psc, hence so does $M'$, by the surgery theorem
of Gromov-Lawson and Schoen-Yau (\cite{GL2} and \cite{SY}---some
of the details are carefully redone in Theorem 3.1 of \cite{RS2}),
and so $Y(M')>0$.  If $M$ is disconnected and 
$Y(M_i)=0$ for some components and $Y(M_j)>0$ for other components,
then we may first replace $M$ by the connected sum of its
components, which has $Y\ge 0$ by iterated application
of case (b) of \cite{PY}, Theorem 1. (See also \cite{Kob}.)
This reduces us to the case where $M$ is connected.
If $M$ is connected and $Y(M)\le0$, then the Corollary to 
Theorem 1 of \cite{PY}
says  $Y(M')\ge Y(M)$, so if $Y(M)$ is exactly $0$, $Y(M')\ge 0$.
\end{demo}
In this paper we will discuss what can be learned about the
\Yi\ for non-simply connected manifolds, using Theorem \ref{thm:surgery}.

Many of the basic facts about manifolds of \psc,
which are proved using the surgery theorem
of Gromov-Lawson and Schoen-Yau, have obvious counterparts for
manifolds with \nYi, obtained by substituting Theorem
\ref{thm:surgery} in the proof.  The proofs are almost identical
to those in the \psc\ case, so while we will give complete statements
of the results, we will be brief when it comes to details of the proofs.

First we need to convert the Surgery Theorem, Theorem
\ref{thm:surgery}, to a Bordism Theorem.  We repeat some definitions
from \cite{RS1} and \cite{RS2}:

\begin{defn} Let $B\to BO$ be a fibration. A \emph{$B$-structure}
on a manifold is defined to be a lifting of the
(classifying map of the) stable normal bundle to a map
into $B$. Then one has bordism groups $\Om^B_n$ of
manifolds with $B$-structures, defined in the usual
way.  (For instance, if $B=B\mathrm{Spin}$,
mapping as usual to $BO$, then
$\Om^B_n=\Oms_n$.) We note that given
a connected closed manifold $M$,
there is a choice of such a $B$
for which $M$ has a $B$-structure and the map $M\to B$
is a 2-equivalence.
\end{defn}
\begin{exs}\label{ex:cases} 
The following special cases show that many of the classical bordism
theories arise via this construction.
\par
\begin{enumerate}
\item If $M$ is a
spin manifold, choose $B=B\pi\times B\mathrm{Spin}$,
where $\pi=\pi_1(M)$, and let $B\to BO$ be the projection onto the second
factor composed with the map 
$B\mathrm{Spin}\to BO$ induced by $\mathrm{Spin}\to O$.
Map $M$ to the first factor
by means of the classifying map for the universal cover,
and to the second factor by means of the spin structure. 
The map $M\to B$ is a $2$-equivalence since it induces
an isomorphism on $\pi_1$ and $\pi_2(B)=0$.
The associated bordism theory is $\Oms_*(B\pi)$.
\item If $M$ is oriented and the universal cover 
$\widetilde M$ of $M$ is non-spin, choose $B=B\pi\times BSO$,
where $\pi=\pi_1(M)$, and let $B=BSO\to BO$ be the obvious map.
Map $M$ to the first factor
by means of the classifying map for the universal cover,
and to the second factor by means of the orientation.
The map $M\to B$ is a $2$-equivalence since it induces
an isomorphism on $\pi_1$ and $\pi_2(B)\cong \pi_2(BSO)
\cong \pi_1(SO)\cong \bZ/2$, with the map $\pi_2(M)\to
\pi_2(B)$ corresponding to $w_2(\widetilde M)$.
The associated bordism theory is $\Om_*(B\pi)$.
\item If $M$ is not orientable and the
universal cover of $M$ is non-spin, let $\pi=\pi_1(M)$, 
and let $B$ be defined by the homotopy pull-back diagram
\[
\xymatrix{B \ar[r] \ar[d] & B\pi \ar[d]^{w_1}\\
BO \ar[r]^{w_1} & \,\bR\bP^\infty,}
\]
where the maps labeled $w_1$ are defined by the first
Stiefel-Whitney class. Note that $BO$ has fundamental group
$\bZ/2$ and that $w_1\co BO\to \bR\bP^\infty$ induces an
isomorphism on $\pi_1$, so that $B$ has fundamental group $\pi$.
The map $B\to BO$ can be taken to be a fibration with fiber
$B\pi'$, where $\pi'=\ker w_1$ is the fundamental
group of the oriented double cover of $M$.
Then the maps of $M$ to $B\pi$
by means of the classifying map for the universal cover and
to $BO$ by means of the classifying map for the
stable normal bundle define a map from $M$ to $B$
which is a $2$-equivalence for the same
reason as in the last example. We will denote the associated
bordism theory by $\fN_*(B\pi\downarrow \bR\bP^\infty)$; it
is a ``twisted version'' of unoriented bordism
with coefficients in $B\pi'$, and it obviously
comes with a natural map to $\fN_*(B\pi)$. In the special case
where $\pi$ splits as $\pi'\times \bZ/2$, with $\pi'=\ker (w_1\co
\pi\to\bZ/2)$, then $B$
becomes simply $B\pi'\times BO$, and the associated bordism theory 
is $\fN_*(B\pi')$. In general,  $\fN_*(B\pi\downarrow \bR\bP^\infty)$
is more complicated to describe, though the following proposition
often tells as much as one needs to know about it.
\end{enumerate}
\end{exs}
\begin{prop}\label{prop:nonorient}In Example \ref{ex:cases}.3,
if $w\co \pi\to\bZ/2$ is surjective,
the natural map $$\fN_*(B\pi\downarrow \bR\bP^\infty)\to \fN_*(B\pi)$$
is injective, at least on classes in dimension $*\ge4$ represented by
disjoint unions of
connected manifolds $M_i$ with non-trivial $w_1$ and with each $u_i\co
M_i\to B\pi$ an isomorphism on $\pi_1$.
\end{prop}
\begin{demo}Consider a
class in $\fN_n(B\pi\downarrow \bR\bP^\infty)$ represented by
$\coprod_i M^n_i$, $n\ge 4$, with $u_i\co M_i^n\to B\pi$ the 
classifying map for the universal cover, and suppose 
the first Stiefel-Whitney class of $M_i$ is $w\circ u_i$ and is
non-trivial. Assume $\coprod_i M_i$ bounds in
$\fN_n(B\pi)$. That means we have a manifold $W^{n+1}$ with $\partial
W^{n+1} = \coprod_i M_i^n$ and with a map $f\co W\to B\pi$ extending 
each $u_i$. We may assume $W$ is connected. Then $f_\sharp$ 
(the induced map on fundamental groups) is a split surjection,
with splitting map $\pi\cong \pi_1(M_i) \stackrel{{\iota_i}_\sharp}
\longrightarrow  \pi_1(W)$, where $\iota_i\co M_i\hookrightarrow W$. 
Since the first Stiefel-Whitney class
for $W$ must extend $w\circ u_i$, it is trivial on $\ker f_\sharp$.
Thus we may do surgery on embedded circles in the interior of $W$ 
and in the kernel of $f_\sharp$ to 
reduce to the case where $f_\sharp$ is also an isomorphism on
$\pi_1$. (The assumption of dimension $*\ge4$ makes it possible
to kill $\ker f_\sharp$ completely with such surgeries.)
Then the first Stiefel-Whitney class of $W$ must be represented by $w\circ f$,
and we obtain a map $W\to B\pi$ extending $M\to B\pi$, showing
that $\coprod (M_i,u_i)$ bounds in $\fN_n(B\pi\downarrow \bR\bP^\infty)$.
\end{demo}

The simply connected cases of the \psc\ analogue of the following theorem
were proved in \cite{GL2}; the general case of the \psc\ analogue, with
this formulation, is in \cite{RS1} and \cite{RS2}.
\begin{thm}[Bordism Theorem]\label{thm:bordthm}
Let $M^n$ be a connected $B$-manifold with $n=\dim M \ge 5$, and assume
that the map $M\to B$ is a 2-equivalence. Then $Y(M)\ge 0$
if and only if the $B$-bordism class of $M$ lies in the
subgroup of $\Om^B_n$ generated by $B$-manifolds with \nYi.
\end{thm}
\begin{sketch} Let $N$ be a $B$-manifold $B$-bordant to $M$. 
The hypotheses combine (via the method of proof of the $s$-Cobordism
Theorem) to show that $M$ can be obtained from $N$ by surgeries
in codimension $\ge 3$. Then if each component of $N$ has \nYi,
one can apply Theorem \ref{thm:surgery} to conclude that the same is true
for $M$. This does it since addition in $\Om^B_n$ comes from connected
sum and additive inverses correspond to reversal of orientation, etc.,
which doesn't affect the \Yi\ of the underlying manifold.
\end{sketch}

Fortunately for applications, one can do better than this.
For simplicity, we restrict attention to the three cases 
discussed in Examples \ref{ex:cases}.
\begin{thm}[Jung, Stolz]\label{thm:homology}
Let $M^n$ be a compact connected manifold with $n=\dim M \ge 5$.
\begin{enumerate}
\item If, as in Example \ref{ex:cases}.1, $M$ is
spin with fundamental group $\pi$, then $Y(M)\ge 0$ if and only
the class  of $M\to B\pi$ in $ko_n(B\pi)$ lies in the subgroup
$ko_n^{\ge 0}(B\pi)$ generated by classes of $M'\to B\pi$
with $M'$ a spin manifold with \nYi, and $M'\to B\pi$ a map 
{\lp}not necessarily an isomorphism on $\pi_1${\rp}. Here $ko_*$ is 
the homology theory corresponding to the connective real $K$-theory
spectrum.
\item If, as in Example \ref{ex:cases}.2, $M$ is oriented
with fundamental group $\pi$, and the universal cover
of $M$
is \textbf{not} spin, then $Y(M)\ge 0$ if and only
the class  of $M\to B\pi$ in $H_n(B\pi,\bZ)$ lies in the subgroup
$H_n^{\ge 0}(B\pi,\bZ)$ generated by classes of $M'\to B\pi$
with $M'$ an oriented manifold with \nYi, and $M'\to B\pi$ a map 
{\lp}not necessarily an isomorphism on $\pi_1${\rp}.
\item If, as in Example \ref{ex:cases}.3, $M$ is non-orientable
with fundamental group $\pi$, and if
the universal cover of $M$
is \textbf{not} spin, then $Y(M)\ge 0$ if and only
the class  of $M\to B\pi$ in $H_n(B\pi,\bZ/2)$ lies in the subgroup
$H_n^{\ge 0}(B\pi,\bZ/2)$ generated by classes of $M'\to B\pi$
with $M'$ a manifold with \nYi, and $M'\to B\pi$ a map.
\end{enumerate}
\end{thm}
\begin{sketch} It was proved by Jung and Stolz (see \cite{RS1} and
\cite{RS2}) that the kernel of the map $\Oms_n(B\pi)\to ko_n(B\pi)$ in case
1, and the kernel of the map $\Om_n(B\pi)\to H_n(B\pi,\bZ)$ in case
2, are represented by manifolds with \psc. Thus the result immediately
follows from Theorem \ref{thm:bordthm}. Now consider Case 3. 
The ``only if'' direction is obvious, so suppose we are given
$M$ non-orientable with  fundamental group $\pi$ and universal cover
non-spin, and assume the class  of $M\to B\pi$ in $H_n(B\pi,\bZ/2)$
lies in the subgroup $H_n^{\ge 0}(B\pi,\bZ/2)$. By Theorem
\ref{thm:bordthm} and Example \ref{ex:cases}.3, it suffices
to show that the class of $M$ in $\fN_n(B\pi\downarrow \bR\bP^\infty)$
lies in $\fN_n^{\ge0}(B\pi\downarrow \bR\bP^\infty)$.
By Proposition \ref{prop:nonorient}, $\fN_n(B\pi\downarrow \bR\bP^\infty)$ is
detected by its image in $\fN_n(B\pi)$. Since
the unoriented bordism spectrum is an Eilenberg-Mac Lane spectrum,
\[
\fN_n(B\pi)\cong \bigoplus_{i+j=n} H_i(B\pi,\bZ/2)\otimes_{\bZ/2} \fN_j,
\]
and each class in $\fN_n(B\pi)$ is a sum of classes of the
form $\left(M_1^i\to B\pi\right) \times M_2^j$, where $M_2^j$
represents a class in $\fN_j$. Here the summand with $j=0$
corresponds to the image of $M\to B\pi$ in $H_n(B\pi,\bZ/2)$. We claim
every class in $\fN_j$, $j>0$, is represented by a manifold
with \nYi. Indeed, multiplicative generators
of $\fN_j$ can be taken to be real projective spaces and 
quadric hypersurfaces in products of  real projective spaces
(\cite{Stong}, p.\ 97).  All of these manifolds admit metrics of \psc\
(cf.\ the argument in the proof of \cite{GL2}, Corollary C),
except for a point in dimension $j=0$. So by Theorem \ref{thm:surgery} above
and Proposition \ref{prop:prod} below, if the class of $M\to B\pi$
lies in $H_n^{\ge 0}(B\pi,\bZ/2)$, then the class of $M\to B\pi$
in $\fN_n(B\pi)$ is represented by a map $M'\to B\pi$,
with $M'$ a manifold with \nYi. Choose
a bordism $f\co W\to B\pi$ between $M\to B\pi$ and $M'\to B\pi$.
As in the proof of Proposition \ref{prop:nonorient}, we may assume
(by doing surgeries on the interior of $W$) that
$f_\sharp$ is an isomorphism on $\pi_1$. 
As in the proof of Proposition \ref{prop:nonorient}, this implies
$M'$ and $M$ represent the same element of 
$\fN_n(B\pi\downarrow \bR\bP^\infty)$, and we conclude using
Theorem \ref{thm:bordthm}.
\end{sketch}

This is now enough machinery to deal with ``easy'' torsion-free
fundamental groups:
\begin{thm} Let $M^n$ be a closed connected $n$-manifold with a fundamental
group $\pi$ which is either free abelian
or of homological dimension $\le 4$. {\lp}This includes
the fundamental groups of  aspherical $2$-manifolds,
$3$-manifolds, and $4$-manifolds.{\rp}
Assume either that $M$ is spin or
that its universal cover is non-spin.  In the spin case, also assume
that the Atiyah-Hirzebruch spectral sequence $H_p(B\pi, ko_q)\Rightarrow
ko_*(B\pi)$ collapses.  {\lp}This is automatic if $\pi$ is
of homological dimension $\le 3$.{\rp}
Then if $n\ge 5$, $M$ has \nYi.
\end{thm}
\begin{demo}
By Theorem \ref{thm:homology}, it's enough to show that for each of
these groups $\pi$, $H_n^{\ge 0}(B\pi,\bZ)$ exhausts $H_n(B\pi,\bZ)$ 
and $ko_n^{\ge 0}(B\pi)$  exhausts $ko_n(B\pi)$ for $n\ge 5$. The
non-spin case is easy, since for $\pi$ free abelian and any $n$,
$H_n(B\pi,\bZ)$ is generated additively by the classes of tori,
which carry flat metrics and thus have \Yi{} zero, whereas if
$\pi$ has homological dimension $\le 4$, $H_n(B\pi,\bZ)$ vanishes for $n\ge 5$.
So consider the spin case.  When $\pi$ is free abelian, the Atiyah-%
Hirzebruch spectral sequence collapses and
\[
ko_n(B\pi)\cong \bigoplus_{p+q=n} H_p(B\pi,\bZ)\otimes ko_q.
\]
Thus this group is generated by the classes of $f\co
T^p\times N^q\to B\pi$, where the map $f$ factors through $T^p$. Since,
as pointed out in \cite{P2},
$ko_*$ is generated by the classes of manifolds of \nYi, we have
the desired result. The other cases are similar
but easier.
\end{demo}

Most of this paper will now deal with the opposite extreme, the case
where $\pi_1(M)$ is finite. In this case, the following results
reduce us to the case where the fundamental group is a $p$-group.
\begin{lemma}\label{lem:coverings}
Suppose $M^n$ is a closed
connected manifold with $Y(M)\ge 0$, and suppose $\widetilde M$
is a finite covering of $M$. Then $Y(\widetilde M)\ge 0$.
\end{lemma}
\begin{demo}Let $m$ be the number of sheets of the
covering $\widetilde M\to M$.
By assumption, given $\varepsilon>0$, we can choose a conformal
class $\cC$ on $M$ with $Y(M,\cC)\ge -\varepsilon$. That
means there is a metric $g$ on $M$ with unit volume and constant
scalar curvature $s\ge -\varepsilon$. Lift the metric $g$ up
to $\widetilde M$. That gives a metric on $\widetilde M$ with
volume $m$ and scalar curvature $s\ge-\varepsilon$. Rescaling, we
get a metric on $\widetilde M$ with unit volume and scalar
curvature $\ge -m^{-\frac 2 n}\varepsilon$. 
This being true for all $\varepsilon>0$,
it follows that  $Y(\widetilde M)\ge 0$.
\end{demo}
\begin{prop}\label{prop:invprop}
If $\pi_1$ and $\pi_2$ are groups and if $\varphi\co \pi_1\to\pi_2$
is a group homomorphism, then $\varphi$ sends $H_n^{\ge 0}(B\pi_1,\bZ)$
to $H_n^{\ge 0}(B\pi_2,\bZ)$, $H_n^{\ge 0}(B\pi_1,\bZ/2)$
to $H_n^{\ge 0}(B\pi_2,\bZ/2)$,  and  $ko_n^{\ge 0}(B\pi_1)$ to
$ko_n^{\ge 0}(B\pi_2)$. If $\pi_1$ is a subgroup of $\pi_2$ of
finite index, then the transfer map on $H_n$ or $ko_n$ sends
$H_n^{\ge 0}(B\pi_2,\bZ)$ to $H_n^{\ge 0}(B\pi_1,\bZ)$,
$H_n^{\ge 0}(B\pi_2,\bZ/2)$ to $H_n^{\ge 0}(B\pi_1,\bZ/2)$, and 
$ko_n^{\ge 0}(B\pi_2)$ to $ko_n^{\ge 0}(B\pi_1)$. 
\end{prop}
\begin{demo}
The first statement is obvious from the definitions in Theorem
\ref{thm:homology}. The second statement follows from
Lemma \ref{lem:coverings}, since the transfer is realized geometrically
via coverings.
\end{demo}
\begin{thm}[Kwasik, Schultz \cite{KS1}]\label{thm:KS}
 Let $M^n$ be a closed connected
$n$-mani\-fold with finite fundamental group $\pi$. 
Assume either that $M$ is spin or
that its universal cover is non-spin.  For each prime $p$, let $i_p\co
\pi_p\hookrightarrow \pi$ be the inclusion of a Sylow $p$-subgroup of $\pi$,
and let $t_p\co H_n(B\pi,\bZ)\to H_n(B\pi_p,\bZ)$, 
$t_p\co H_n(B\pi,\bZ/2)\to H_n(B\pi_p,\bZ/2)$, $t_p\co ko_n(B\pi)\to 
ko_n(B\pi_p)$ be the transfer maps. Then $M$ has \nYi\ if and only if
$t_p([M])$ lies in the subgroup $H_n^{\ge 0}(B\pi_p,\bZ)$ in the
oriented non-spin case, or in
$ko_n^{\ge 0}(B\pi)$ in the spin case, for each $p$ dividing 
the order of $\pi$. In the non-orientable  non-spin case, 
$M$ has \nYi\ if and only if
$t_2([M])$ lies in $H_n^{\ge 0}(B\pi_2,\bZ/2)$.
\end{thm}
\begin{demo}The proof is almost word-for-word as in \cite{KS1},
but we review the argument. The ``only if'' statement
is contained in Proposition \ref{prop:invprop}. As for the
``if'' statement, let $A=\widetilde H_n(B\pi,\bZ)$, 
$\widetilde H_n(B\pi,\bZ/2)$, or $\widetilde 
ko_n(B\pi)$, and let $B$ be the subgroup $\widetilde H_n^{\ge 0}
(B\pi,\bZ)$,  $\widetilde H_n^{\ge 0}(B\pi,\bZ/2)$, or $\widetilde 
ko_n^{\ge 0}(B\pi)$. Similarly let $A_p=\widetilde H_n(B\pi_p,\bZ)$,
$\widetilde H_n(B\pi_p,\bZ/2)$, or $\widetilde 
ko_n(B\pi_p)$, and let $B_p$ be the subgroup $\widetilde H_n^{\ge 0}
(B\pi_p,\bZ)$, $\widetilde H_n^{\ge 0}(B\pi_p,\bZ/2)$, 
or $\widetilde ko_n^{\ge 0}(B\pi_p)$. 
(We can work with reduced homology since
$H_*(\text{pt})=H_*^{\ge 0}(\text{pt})$ and $ko_*(\text{pt})=
ko_*^{\ge 0}(\text{pt})$.) 

Note that $A$ is a finite group and $B$
is a subgroup; we are trying to show that an element $[M]$ of $A$
lies in $B$ if $t_p([M])\in B_p$ for all $p$. Now
$\alpha_p=i_p\circ t_p$ is an endomorphism of $A$ which is an
isomorphism on $A_{(p)}$, since $[\pi:\pi_p]$ is a unit modulo $p$.
If $t_p([M])\in B_p$ for all $p$, then $\alpha_p([M])\in i_p(B_p)\subseteq
B$ for all $p$, by Proposition \ref{prop:invprop}. So that means the
image of $[M]$ in $A_{(p)}$ lies in $B_{(p)}$ for all $p$,
and thus $[M]$ lies in $B$. Note, incidentally, that in the
non-orientable, non-spin case, only the Sylow 2-subgroup matters,
since $\widetilde H_*(B\pi_p,\bZ/2)=0$ for $p\ne 2$.
\end{demo}
\begin{thm}[Kwasik, Schultz \cite{KS1}]\label{thm:stablesplitting} 
Let $\pi$ be a finite group, and let 
\[
e\co \Omega^\infty\Sigma^\infty B\pi_+
\to \Omega^\infty\Sigma^\infty B\pi_+
\]
be an idempotent in the stable homotopy category
{\lp}giving a stable splitting of $B\pi${\rp}. Then $e$ maps
$H_n^{\ge 0}(B\pi,\bZ)$,
$H_n^{\ge 0}(B\pi,\bZ/2)$, and $ko_n^{\ge 0}(B\pi)$ into themselves.
\end{thm}
\begin{sketch}
As pointed out in \cite{KS1}, the proof of the Segal Conjecture
implies that the stable splittings of $B\pi$ are essentially
linear combinations of products of transfer maps and maps induced by
group homomorphisms, so the result then follows from
Proposition \ref{prop:invprop}. 
\end{sketch}

\section{Analytic Tools}\label{sec:tools}
In this section we present a number of analytic results that can be
used to study the classes of manifolds with \psc\ or with \nYi.
First we need a basic characterization of manifolds in the latter
class.
\begin{prop}\label{prop:nnYchar}Let $M^n$ be a closed $n$-manifold.
Then:

{\lp i\rp} If $M$ does not admit a metric of \psc, then
\begin{equation}\label{eq:Ycalc}
Y(M)=-\inf_g \left(\int_M |s_g|^{\frac n 2}\dvol_g\right)^{\frac 2 n}.
\end{equation}
Here the infimum is taken over all Riemannian metric $g$ on $M$,
and $s_g$ denotes the \scurv\ of $g$.

{\lp ii\rp} Suppose that for each $\varepsilon>0$, there exists a 
metric $g$ on $M$ with volume $1$ and $|s_g|<\varepsilon$. Then
$Y(M)\ge 0$. The converse is true if $n\ge 3$ or if $Y(M)= 0$.
\end{prop}
\begin{demo}Statement (i) is Proposition 1 in \cite{LeB}.
As for (ii), suppose the condition is satisfied. If $Y(M)>0$, then
we have nothing to prove, and if not, (i) shows that $Y(M)\ge 0$.
In the converse direction, suppose  $Y(M)\ge 0$. If $n\ge 3$ and 
if $Y(M)>0$, then by a theorem of Kazdan and Warner \cite{KW},
$M$ admits a metric $g$ with $s_g\equiv 0$, and obviously we may
rescale $g$ to have volume 1 without changing this condition. 
If, on the other hand, $Y(M)=0$, that means, by
definition of the \Yi\ (recall
equation \eqref{eq:Ydef}), that for all $\varepsilon>0$, there exists a 
metric $g$ on $M$ with volume $1$ and $s_g\le 0$ constant and $>
-\varepsilon$. So again the condition of (ii) is satisfied.
\end{demo}

Another basic fact is the following:
\begin{prop}\label{prop:prod}
Suppose $M^m$ and $N^n$ are closed manifolds, $n=\dim N\ge 1$, and
also $Y(N)\ge 0$. Then $Y(M\times N)\ge 0$.
\end{prop}
\begin{demo}If $Y(N)>0$, then $N$ admits a metric of \psc\ and so does
$M\times N$, so $Y(M\times N)> 0$. If $Y(N)=0$, then by Proposition
\ref{prop:nnYchar}, given $\varepsilon>0$, there exists a 
metric $g_\varepsilon$ on $N$ with volume $1$ and $|s_g|<\varepsilon$. Choose
any metric $g'$ on $M$ with volume $1$. If we give $M\times N$
the product metric $g\times tg'$ (where $tg'$ means $g'$ rescaled by
multiplying all distances by $t$), then this product metric has 
\scurv\ $t^{-2}s_{g'}+s_{g_\varepsilon}$ and volume $t^m$. So the
integral in \eqref{eq:Ycalc} (with $M\times N$ in place of $M$)
is
\begin{equation}\label{eq:prodest}
\iint \left\vert t^{-2}s_{g'}+s_{g_\varepsilon}\right\vert^{\frac{n+m}{2}}
\dvol_{tg'}\dvol_\varepsilon 
\le t^m \left\vert C
t^{-2}+\varepsilon\right\vert^{\frac{n+m}{2}},
\end{equation}
for some constant $C$ (the maximum of $|s_{g'}|$ over $M$)
independent of $t$ and $\varepsilon$. So the idea is to take
$t$ large, and then given $t$, to take $\varepsilon$ of order
$t^{-2}$. In equation \eqref{eq:prodest}, we see that
the integral on the left-hand side is bounded by a constant times
\[t^m \left( t^{-2} \right)^{\frac{n+m}{2}} = t^{-n},
\] 
which goes to zero as $t\to \infty$. Thus by Proposition
\ref{prop:nnYchar}, the result follows.
\end{demo}

Next, we discuss the extension of the minimal hypersurface technique
of \cite{SY} to the study of \nYi. Suppose $M^n$ is a closed manifold
with metric $g$, and suppose $H^{n-1}$ is a stable minimal hypersurface
in $M$. In \cite{SY}, it was shown that if $s_g>0$, then 
$Y(H, [\bar g])>0$, where $\bar g$ denotes the induced metric on $H$
and $[\bar g]$ is its conformal class. In particular, there is a
metric in  $[\bar g]$ with \psc, and this can be used to rule out
\psc\ metrics on many non-simply connected manifolds.
Now it is not true that just because $Y(M)\ge 0$, then
$Y(H)\ge 0$, since by Proposition \ref{prop:prod},
we can get a counterexample by taking $M=H\times S^1$
and $Y(H)<0$ (say with $n-1=2$ or $4$). However, the same estimates
used in the proof of Theorem in \cite{SY} show that if
$s_g\ge K$, where $K$ is a constant, then
because the second variation of the $(n-1)$-dimensional
volume of $H$ is positive, one has
\begin{equation}\label{eq:secvar}
\int_H \frac{(\bar s - K)\phi^2}2 + \int_H|\nabla\phi|^2 >0  
\end{equation}
for all functions $\phi\in C^\infty(H)$ not vanishing identically.
(Here $\bar s$, the scalar curvature of $H$, and $\nabla$
are to be computed with respect to the induced metric $\bar g$.) 
Assume $n>3$ and consider the ``conformal Laplacian''
\[
L_H =\frac{4(n-2)}{n-3}\Delta_H  + \bar s
\]
of $H$, where $\Delta_H$ is the usual (non-negative) Laplacian. 
(Recall that the dimension of $H$ is $n-1$, not $n$.) Then
for $\phi$ as above we have
\begin{align}
\frac12 \langle L_H\phi,\,\phi\rangle&=
      \frac{2(n-2)}{n-3}\int_H|\nabla\phi|^2
      +\int_H \frac{\bar s\phi^2}2 \notag\\
      &=\frac{n-1}{n-3}\int_H|\nabla\phi|^2+
      \left(\int_H \frac{(\bar s - K)\phi^2}2 + \int_H|\nabla\phi|^2
      \right) +\frac{K}{2}\int_H \phi^2 \notag\\
      &> \frac{K}{2}\int_H \phi^2.
\end{align}
Note the use of equation \eqref{eq:secvar} at the last step.
This implies that if $K$ is close to $0$, then the conformal
Laplacian $L_H$ is not too negative, and thus $Y(H, \bar g)$
is not too negative, provided that the ($n-1$)-dimensional
volume of $H$ is not too large.
                                                             
If $n=3$, things are even easier: we instead take $\phi\equiv1$ in 
equation \eqref{eq:secvar} and apply Gauss-Bonnet.
These arguments thus prove the following two results:
\begin{thm}\label{thm:minhyp} Let $M^n$ be a closed manifold
with metric $g$, and suppose $H^{n-1}$ is a stable minimal hypersurface
in $M$. Also suppose that the metric $g$ is scalar-flat. Then
$Y(H)\ge 0$.
\end{thm}
\begin{demo}Immediate from the above estimates.
\end{demo}
\begin{thm}\label{thm:circcrosssurf}
Let $M^2$ be a closed oriented surface of genus $g>1$, and let
$N^3=S^1\times M^2$. Then $Y(N)=0$ by Proposition \ref{prop:prod}.
{\lp}It cannot be strictly positive, by \textup{\cite{GL3}}, Theorem
8.1, for example.{\rp} Thus there is a sequence $g_i$ of metrics
on $N^3$ with volume $1$ and constant scalar curvatures
$s_i$, with the \scurv s tending to $0$ as $i\to \infty$. On the
other hand, for any such sequence of metrics, $\diam(N, g_i)\to
\infty$.
\end{thm}
\begin{demo}Choose the metrics $g_i$ as in the statement of
the theorem.  Choose minimal submanifolds $M_i$
which are absolutely area-minimizing in the homology class
$[M^2]\in H_2(S^1\times M^2,\bZ)$ for the metric $g_i$. 
By inequality \eqref{eq:secvar} with $\phi\equiv 1$, 
\begin{equation}\label{eq:GBlimit}
\liminf_{i\to\infty}
\int_{M_i} \left(\bar s_i- s_i\right)\dvol_{\bar g_i} \ge 0.
\end{equation}
On the other hand, each $M_i$ must be a surface of genus $>1$,
since it represents a non-trivial homology class in the
infinite cyclic cover $\bR\times M^2$, while each mapping
of a sphere into this space is null-homotopic and each
mapping of a torus into this space factors through a circle
(since each abelian subgroup of $\pi_1(M)$ is cyclic)
and is thus trivial in $H_2$.
So by Gauss-Bonnet, $\int_{M_i} \bar s_i\dvol_{\bar g_i}
\le -4\pi$. Comparing this with equation \eqref{eq:GBlimit},
we see the area of $M_i$ with respect to $\bar g_i$ must tend
to $\infty$ as $i\to\infty$, while the average value
of $\bar s_i$ must go to $0$, and in particular, $\diam (M_i,
\bar g_i)\to\infty$. This in turn means $\diam g_i\to\infty$,
since otherwise we could choose representatives for the homology
class $[M]$ in $(N, g_i)$ with bounded diameters, a contradiction.
\end{demo}

The next two results are the most significant
in this paper; they will be used in the next section to deal
with ``Toda brackets,'' among the most intractable of bordism classes.

\begin{thm}\label{thm:pscToda}Let $M_0$ and $M_1$ be closed manifolds,
not necessarily connected, that admit metrics of \psc.  Suppose
$M_0=\partial W_0$ and $M_1=\partial W_1$ for some compact manifolds
with boundary, $W_0$ and $W_1$. Form a new manifold
\[ M = \left(W_0\times M_1\right) \cup_{M_0\times M_1}
\left(M_0\times W_1\right) 
\]
of dimension $n_0+n_1+1$, where $n_0$ and $n_1$ are the dimensions of
$M_0$ and $M_1$. Then $M$ admits a metric of \psc.
\end{thm}
\begin{demo}We start by choosing metrics of \psc, $g_0$ and $g_1$,
on $M_0$ and $M_1$, respectively. Extend them to metrics $\bar g_0$
and $\bar g_1$ on $W_0$ and $W_1$, which are product
metrics in neighborhoods of the boundaries. The trick is to write $M$ as a 
union of four pieces (not two) as follows:
\[ M = \left(W_0\times M_1\right) \cup_{M_0\times M_1}
\left(T_0\times M_1\right) \cup_{M_0\times M_1}
\left(M_0\times T_1\right)\cup_{M_0\times M_1}
\left(M_0\times W_1\right),
\]
where the ``tubes'' $T_0$ and $T_1$ are (as smooth manifolds)
$M_0\times I$ and $M_1\times I$, respectively.
Call the pieces here $A_0$, $T_0\times M_1$, $M_0\times T_1$, and 
$A_1$ in that order. Since $g_0$ and $g_1$ have \psc, we can
choose (very small) constants $t_0>0$ and $t_1>0$ so that
the metric $\bar g_0 \times t_1g_1$ on $A_0$ and 
the metric $t_0g_0\times \bar g_1 $ on $A_1$ have \psc. Now all
we have to do is choose the metric $g_{T_0}$ on $T_0$ to interpolate between
$t_0 g_0$ and $g_0$ and the metric $g_{T_1}$ on $T_1$ to interpolate between
$t_1 g_1$ and $g_1$. If the tubes $T_0$ and $T_1$ are ``very long,''
it is possible to do this so that  $T_0$ and $T_1$ have \psc,
by the ``Isotopy implies concordance'' theorem, \cite{GL2}, Lemma 3.
(In fact, in this case, one can write down an explicit warped product
metric that does the trick.)
Then all the metrics fit together to give a metric of \psc\ on $M$.
\end{demo}

The next theorem is quite similar, but considerably more delicate.

\begin{thm}\label{thm:nnYiToda}Let $M_0$ and $M_1$ be closed manifolds,
not necessarily connected, each with \nYi.  {\lp}When
$M_i$ is disconnected, we mean that each component is required
to have \nYi.{\rp} Suppose
$M_0=\partial W_0$ and $M_1=\partial W_1$ for some compact manifolds
with boundary, $W_0$ and $W_1$. Form a new manifold
\[ M = \left(W_0\times M_1\right) \cup_{M_0\times M_1}
\left(M_0\times W_1\right) 
\]
of dimension $n_0+n_1+1$, where $n_0$ and $n_1$ are the dimensions of
$M_0$ and $M_1$. Then, excluding the case where $Y(M_0)=0$,
$n_1=2$, and $Y(M_1)>0$, it follows that $Y(M)\ge 0$.
\end{thm}
\begin{demo} We follow the same approach as in the proof of
Theorem \ref{thm:pscToda}. If $Y(M_0)$ and $Y(M_1)$ are both
strictly positive, we're done by Theorem \ref{thm:pscToda},
so we may assume at least one of $M_0$ and $M_1$ has $Y=0$.
Then we're excluding the case where the other manifold is 
a disjoint union of copies of $S^2$
or $\bR\bP^2$, so by Proposition \ref{prop:nnYchar}, we may assume
both manifolds have metrics of unit volume which are almost
scalar-flat.
By Proposition \ref{prop:nnYchar}, it is 
enough to show that $M$ has a metric for which the integral in
\eqref{eq:Ycalc} is as small as one likes. We will estimate the
integral separately over the four pieces of $M$ (as in the
last proof) and add the results. Fix $\varepsilon>0$
and choose metric $g_0$ and $g_1$
on $M_0$ and $M_1$, respectively, each with volume $1$ and 
with small constant scalar curvatures, $s_0$ and $s_1$, respectively,
with $|s_0|, \, |s_1|<\varepsilon$.
Extend $g_0$ and $g_1$ to metrics $\bar g_0$
and $\bar g_1$ on $W_0$ and $W_1$, which are product
metrics in neighborhoods of the boundaries. Then the scalar curvature
of the metric $\bar g_0\times t_1g_1$
on $A_0$ is $s_{\bar g_0}+t_1^{-2}s_1$ and the scalar curvature
of the metric $t_0g_0 \times \bar g_1$
on $A_1$ is $s_{\bar g_1}+t_0^{-2}s_0$. (The constants $t_0$
and $t_1$ will be chosen later.) Furthermore, the volumes of these
metrics are $\vol(\bar g_0)\times t_1^{n_1}$ for $A_0$ and
$\vol(\bar g_1)\times t_0^{n_0}$ for $A_1$. Letting $t_0$
and $t_1$  go to $0$, we see there are constants $c_0>0$ and
$c_1>0$ with 
\begin{align}
\iint\limits_{A_0} \left\vert s_{\bar g_0}+t_1^{-2}s_1
\right\vert^{\frac{n_0+n_1+1}2}
\dvol_{\bar g_0}\dvol_{t_1g_1}
&\le c_0 t_1^{-(n_0+n_1+1)}\varepsilon^{\frac{n_0+n_1+1}2}t_1^{n_1},
\label{eq:A0int}\\
\iint\limits_{A_1} \left\vert s_{\bar g_1}+t_0^{-2}s_0
\right\vert^{\frac{n_0+n_1+1}2}
\dvol_{t_0 g_0}\dvol_{\bar g_1}
&\le c_1 t_0^{-(n_0+n_1+1)}\varepsilon^{\frac{n_0+n_1+1}2}t_0^{n_0}.
\label{eq:A1int}
\end{align}
The right-hand sides of \eqref{eq:A0int} and
\eqref{eq:A1int} can be rewritten as
\[
c_0 t_1^{-(n_0+1)}\varepsilon^{\frac{n_0+n_1+1}2}
=
c_0 \varepsilon^{n_1/2}\left(\frac\varepsilon{t_1^2}\right)^{\frac{n_0+1}2}
\]
and
\[
c_1 t_0^{-(n_1+1)}\varepsilon^{\frac{n_0+n_1+1}2}
=
c_1 \varepsilon^{n_0/2}\left(\frac\varepsilon{t_0^2}\right)^{\frac{n_1+1}2},
\]
respectively.

Next we need to deal with the tubes $T_0$ and $T_1$. We give these
warped product metrics of the form $f_i(x)g_i\times g_\bR$, $i=0,\,1$, where
$g_\bR$ is the standard metric on the line, and $x$ is the parameter
along the length of the tube. The function $f_i$ will be chosen
to interpolate between $0$ and $t_i$. If we write $f_i=\exp(-u_i)$,
we need to choose $u_i$ as in the following picture, so that
the graph has vanishing first and second derivatives at both ends:
\begin{center}\epsfig{file=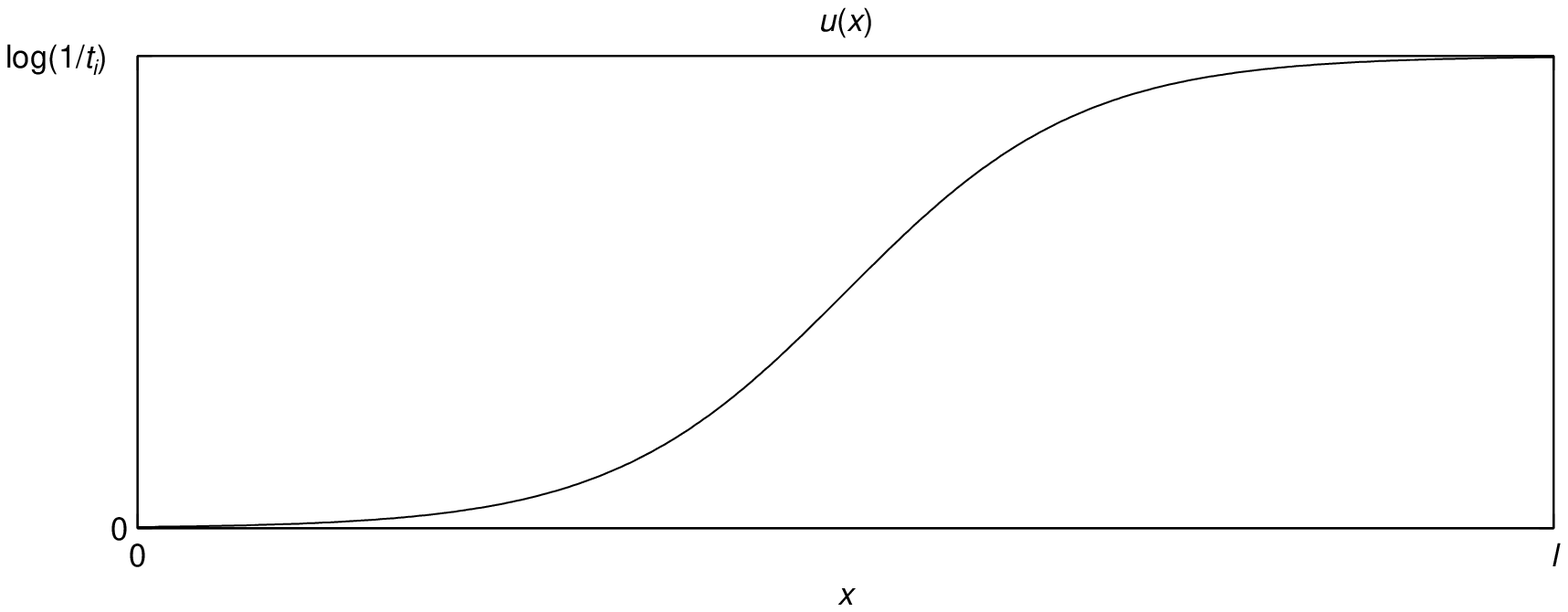,width=300pt}\end{center}
Here $l$, to be taken large, is the length of the tube.
Since $t_i<1$ and $\vol(g_i)=1$, 
the volume of $T_i$ will be bounded by $l$, as will the volume
of $T_0\times M_1$ or $M_0\times T_1$, when we take the product
with the metric $g_1$ on $M_1$ or $g_0$ on $M_0$. The scalar
curvature of $T_i$ is given by equation (7.35) on p.\ 157 of \cite{GL3}, which
gives:
\[
\frac{1}{f_i^2}s_i - \frac{n_i(n_i-1)}{f_i^2}({f_i}')^2-\frac{2n_i}{f_i}{f_i}''.
\]
Since $f=\exp(-u_i)$, $\frac{{f_i}'}{f_i}=-{u_i}'$ and
$\frac{{f_i}''}{f_i}=\left({u_i}'\right)^2-{u_i}''$. Now we can choose
$u_i$ so that ${u_i}'$ is bounded by a constant times
$\frac{\log(1/t_i)}{l}$ and ${u_i}''$ is bounded by a constant times
$\frac{\log(1/t_i)}{l^2}$. Thus the scalar curvature of
$T_i$ is bounded in absolute value by
\[\frac{\varepsilon}{t_i^2} + d_i \frac{\left(\log(1/t_i)\right)^2}{l^2}
\]
for some constant $d_i$. Thus the integrals over $T_0\times M_1$ and 
$M_0\times T_1$ give:
\begin{align}
\iint\limits_{T_0\times M_1} \left\vert s_{T_0}+s_1
\right\vert^{\frac{n_0+n_1+1}2}
\dvol_{T_0}\dvol_{g_1}
&\le l \left\vert \varepsilon + \frac{\varepsilon}{t_0^2} + 
d_0 \frac{\left(\log(1/t_0)\right)^2}{l^2}\right\vert^{\frac{n_0+n_1+1}2}
\label{eq:T0int}\\
\iint\limits_{M_0\times T_1} \left\vert s_{g_0}+s_{T_1}
\right\vert^{\frac{n_0+n_1+1}2}
\dvol_{g_0}\dvol_{T_1}
&\le l \left\vert \varepsilon + \frac{\varepsilon}{t_1^2} + 
d_1 \frac{\left(\log(1/t_1)\right)^2}{l^2}\right\vert^{\frac{n_0+n_1+1}2}.
\label{eq:T1int}
\end{align}
Now all we have to do is choose the parameters $t_0$, $\varepsilon$,
and $l$ to make all of \eqref{eq:A0int}, \eqref{eq:A1int}, 
\eqref{eq:T0int}, and \eqref{eq:T1int} simultaneously small.  We do
this as follows.  First choose $t_0$ and $t_1$ very small.  Then
choose $l$ large enough so that the terms
\[ l \left\vert \frac{\left(\log(1/t_i)\right)^2}{l^2}
\right\vert^{\frac{n_0+n_1+1}2} =
\frac {\left(\log(1/t_i)\right)^{n_0+n_1+1}}{l^{n_0+n_1}}
\]
are small.  Then finally choose $\varepsilon/t_i^2$ extremely small so that 
\[
l \left(\frac{\varepsilon}{t_i^2}\right)^{\frac{n_0+n_1+1}2}
\]
is also small.  That does it.
\end{demo}

\section{Applications to
Non-Negativity of the \\ Yamabe Invariant}\label{sec:Todaresults}
We're now ready for the first main results of this paper.
\begin{thm}\label{thm:nnYinonspinabel}
Let $M^n$ be a closed, connected, $n$-manifold with abelian
fundamental group, with non-spin universal cover, and with $n\ge 5$.
Then $M$ has \nYi.
\end{thm}
\begin{demo}Consider the oriented case. (The 
non-orientable case works exactly the same way, but with integral
homology replaced by homology with coefficients in $\bZ/2$.)
By Theorem \ref{thm:homology}, it's enough to show that 
$H_*^{\ge 0}(B\pi_1(M),\bZ)$ exhausts the image in
$H_*(B\pi_1(M),\bZ)$ of $\Om_*(B\pi_1(M))$.
Write $\pi_1(M)$ as $\bZ^k\times \pi$, with $\pi$
finite abelian. Since the homology of a free abelian group is
torsion free, the K\"unneth Theorem gives
\[
H_n(B\pi_1(M),\bZ)\cong \bigoplus_{p+q=n}H_p(B\bZ^k,\bZ)\otimes
H_q(B\pi,\bZ),
\]
and so the homology of $B\pi_1(M)$ is generated by classes of
products of tori with homology classes of $B\pi$. So by
Proposition \ref{prop:prod}, we only have to show that
$H_*^{\ge 0}(B\pi,\bZ)$ exhausts the image in $H_*(B\pi,\bZ)$
of $\Om_*(B\pi_1(M))$. In other words,
we are reduced to the case of \emph{finite} abelian groups.
By Theorem \ref{thm:KS}, we can further assume that $\pi_1(M)$ 
is a finite abelian $p$-group for some prime $p$ (and in the
non-orientable case, we can further assume that $p=2$).
We will come back to finite abelian $p$-groups after a short digression.
\end{demo}
\begin{lemma}\label{lem:cyclic}Let $\pi$ be a cyclic group
of prime power order $p^k$. Then each class in
$H_n(B\pi,\bZ)$ is represented by an oriented manifold with \nYi, and
if $n>1$, by an oriented manifold with \psc. If $p=2$, then also each class in
$H_n(B\pi,\bZ/2)$ is represented by a manifold 
{\lp}not necessarily orientable{\rp} with \nYi, and
if $n>1$, by a manifold with \psc.
\end{lemma}
\begin{demo}Note that $H_{2n+1}(B\pi,\bZ)$ is cyclic of order $p^k$,
with a generator represented by the lens space $S^{2n+1}/\pi\to
B\pi$, and  $H_{2n}(B\pi,\bZ)$ vanishes for $n>0$. Since the lens space has
{\psc} except in the exceptional case $n=0$, when it has a flat metric,
the first statement is immediate.

Now consider the case of $\bZ/2$-homology and $p=2$. Then
$H_n(B\pi,\bZ/2)$ is cyclic of order $2$ for all $n\ge 1$. When $k=1$,
things are again easy, as we have generators $\bR\bP^n\to B\pi$
for all the homology groups, with {\psc} except in the exceptional
case $n=1$. This leaves the case $k>1$. In that case, 
$H^*(B\pi,\bZ/2)\cong \bF_2[u,v]/(v^2)$, where $\bF_2=\bZ/2$ is
the field of two elements, $v\in H^1(B\pi,\bF_2)$ corresponds to
the group homomorphism $\pi\twoheadrightarrow \bZ/2$, and
$u\in H^2(B\pi,\bF_2)$ is the class of the central extension
\[
1\to \langle a^{2^k}\mid  a^{2^{k+1}}=1\rangle \to
\langle a\mid  a^{2^{k+1}}=1\rangle \to
\langle a\mid  a^{2^{k}}=1\rangle \to 1,
\]
via the usual correspondence between $H^2$ and central extensions.
Since the pull-back of this central extension to the 
unique two-element subgroup of $\pi$ is non-trivial, the inclusion
$\iota\co \bZ/2\hookrightarrow \pi$ induces an isomorphism 
\[
\iota^*\co H^2(B\pi,\bF_2)\to H^2(B\bZ/2,\bF_2),
\]
and thus $\iota^*(u^n)$ is the generator of $H^{2n}(B\bZ/2,\bF_2)$.
By duality, 
\[
\iota_*\co H_{2n}(B\bZ/2,\bF_2)\to H_{2n}(B\pi,\bF_2)
\]
is an isomorphism, and so the generator of $H_{2n}(B\pi,\bF_2)$
is represented by $$\bR\bP^{2n}\to B\bZ/2 \stackrel{B\iota}
{\longrightarrow} B\pi, \ \ \ n\ge1.$$ We still need to find geometric
generators for $H_{2n+1}(B\pi,\bF_2)$. Since $H^{2n+1}(B\pi,\bF_2)$
is generated by $u^nv$, and the cup product comes from restricting
the exterior product $u^n\boxtimes v\in H^{2n+1}(B\pi\times
B\pi,\bF_2)$ to the diagonal copy of $B\pi$, we see by duality
that the generator of $H_{2n+1}(B\pi,\bF_2)$ is represented by
$\Delta^*(\bR\bP^{2n}\times S^1\to B\pi\times B\pi)$, where
$\Delta\co B\pi\to B\pi\times B\pi$ is the diagonal map and
$\Delta^*$ denotes the transfer (which corresponds to taking a 
covering). Since transfer preserves the {\psc} or the {\nYi} property
by Proposition \ref{prop:invprop}, the proof is complete.
\end{demo}
\begin{labeleddemo}{Proof of Theorem \ref{thm:nnYinonspinabel}, continued}
Recall that we have already reduced to the case where the fundamental
group $\pi$ of $M$ is a finite abelian $p$-group, hence a finite product of
cyclic $p$-groups. In the non-orientable case, since $H_*(B\pi,\bF_2)$
is by the K\"unneth Theorem generated by products of homology classes
for cyclic groups, and since all these homology classes are
represented by manifolds of {\nYi} by Lemma \ref{lem:cyclic}, we are done
by Proposition \ref{prop:prod}. In the orientable case, things are
more complicated because we have to deal with the Tor terms in
the K\"unneth Theorem, and also because the natural map (the Hurewicz
homomorphism for $MSO$) $\Om_*(B\pi)\to H_*(B\pi,\bZ)$ may
not be surjective, and may not be split onto its image.
Thus the argument will require some care. 
We prove the theorem by induction on the rank
(the number of cyclic factors in a product decomposition) of $\pi$. If
the rank is $1$, $\pi$ is cyclic and we are done by Lemma
\ref{lem:cyclic}. So assume the result is true for $p$-groups of
smaller rank, and write $\pi=\pi'\times \bZ/p^k$, where we may assume
that $p^k$ is less than or equal to the order of every cyclic factor
of $\pi'$, and thus less than or equal to the order of every cyclic
factor of the homology of $B\pi'$. First assume that $p=2$. This
case is easier because $MSO$ localized at $2$ is a direct sum of
Eilenberg-Mac Lane spectra (see \cite{MM} and \cite{Tay},
or \cite{RosCob} for a review of the literature),
and thus $H_*(B\pi,\bZ)$ can be identified with a direct summand in
$\Om_*(B\pi)$, and similarly for $\pi'$. By inductive hypothesis,
each cyclic factor (say of order $p^s$, $s\ge k$)
in $H_j(B\pi',\bZ)$ is generated by the 
class of a manifold $M'\to B\pi'$, where $M'$ has {\nYi}. If $n-j$ is
odd, then we get a corresponding tensor term in the K\"unneth formula
for the homology of $B\pi$, and it is represented by a product of $M'$
with either $S^1$ or a lens space, and so it is represented by a
manifold with {\nYi}. If $n-j\ge 2$ is even, there is a contribution
to $H_n(B\pi,\bZ)$ of the form $\Tor_\bZ(\bZ/p^s, \bZ/p^k)$, coming from
$H_j(B\pi',\bZ)$ and $H_{n-j-1}(B\bZ/p^k,\bZ)$, which we need
to represent by a manifold of {\nYi}. Since Tor is left exact and
$\Tor_\bZ(\bZ/p^s, \bZ/p^k)$ is cyclic of order $p^k$, the map
\[
\Tor_\bZ(\bZ/p^k, \bZ/p^k) \to \Tor_\bZ(\bZ/p^s, \bZ/p^k)
\]
induced by the inclusion $\bZ/p^k \hookrightarrow \bZ/p^s$ is an
isomorphism, so without loss of generality, we may replace $M'$ by
something representing a multiple of its homology class, and assume
$s=k$. Choose $M''\to B\bZ/p^k$, with
$M''$ either $S^1$ or a lens space, of
dimension ${n-j-1}$, generating $H_{n-j-1}(B\bZ/p^k,\bZ)$. We may
assume the bordism classes of $M'\to B\pi'$ and $M'' \to B\bZ/p^k$ both
have order $p^k$. Then their Tor product in the homology of $B\pi$ may
be represented by the cobordism Massey product $\langle M', p^k, M''
\rangle$ (see \cite{Al}), or in other words, by a Toda bracket 
construction as in
Theorem \ref{thm:nnYiToda}. More precisely,
choose $W_0$ bounding $p^k M'$ over
$B\pi'$ and $W_1$ bounding $p^k M''$ over
$B\bZ/p^k$, and glue together $W_0\times M''$ and $M'\times W_1$ along
their common boundary. By Theorem \ref{thm:nnYiToda}, the resulting
manifold is represented by a manifold with {\nYi}. (Note that the
exceptional case of that theorem never arises.) This completes the
inductive step.

Now consider the case where $p$ is odd. In this case, it's important
to point out that the inductive hypothesis is simply that the \emph{image}
of $\Om_*(B\pi')\to H_*(B\pi',\bZ)$ is represented by manifolds with
{\nYi}, as this map is not usually surjective. However, in this
case we have one additional tool in our arsenal, namely Landweber's
K\"unneth Theorem for oriented bordism \cite{Landw}. More
precisely, we apply Theorem A of \cite{Landw}, which applies since
$\widetilde H_*(B\bZ/p^k, \bZ)$ consists entirely of odd torsion and the
Atiyah-Hirzebruch spectral sequence $\widetilde H_*(B\bZ/p^k,\Om_*)
\Rightarrow \widetilde \Om_*(B\bZ/p^k)$ collapses for dimensional reasons.
($\Om_*$ localized at $p$ is free over $\bZ_{(p)}$
and concentrated in degrees divisible by $4$, and  $\widetilde
H_*\bigl(B\bZ/p^k, \bZ_{(p)}\bigr)$ is non-zero only in odd degrees.)
Note also that the proof of Landweber's Theorem shows that
$\Om_*(B\bZ/p^k)_{(p)}$ has homological dimension $1$ over
$\Om_*(\text{pt})_{(p)}$. Now observe that we have a commutative
diagram with exact rows:
\begin{equation}\label{eq:Torchase}
\def\objectstyle{\scriptstyle}
\xymatrix@C-9pt{
0\ar[r] \ar@{=}[d] & \Om_*(B\pi')\otimes_{\Om_*} \Om_*(B\bZ/p^k)
\ar[r] \ar[d]^\alpha & \Om_*(B\pi) \ar[r] \ar[d]^\beta & 
\Tor_{\Om_*}\bigl(\Om_*(B\pi'),\Om_*(B\bZ/p^k)\bigr) \ar[r] 
\ar[d]^\gamma &
0 \ar@{=}[d] \\
0\ar[r] & H_*(B\pi')\otimes_{\bZ} H_*(B\bZ/p^k)
\ar[r]  & H_*(B\pi) \ar[r]  & 
\Tor_{\bZ}\bigl(H_*(B\pi'),H_*(B\bZ/p^k)\bigr) \ar[r]&
\,0,
}
\end{equation}
in which the vertical arrows are induced by the natural transformation
$\Om_* \to H_*$. Note that the map $\Om_*(B\bZ/p^k)\to H_*(B\bZ/p^k)$
is surjective, and denote the image of $\Om_*(B\pi')\to H_*(B\pi')$ by
$RH_*(B\pi')$ (for ``representable homology'').  The image of $\alpha$
is then obviously $RH_*(B\pi')\otimes_{\bZ} H_*(B\bZ/p^k)$; classes
here are represented by products of manifolds of {\nYi} (because of
the inductive hypothesis), so these have {\nYi}.  The image of $\beta$
is by definition $RH_*(B\pi)$, whereas the image of $\gamma$ is
evidently contained in
$\Tor_{\bZ}\bigl(RH_*(B\pi'),H_*(B\bZ/p^k)\bigr)$. Also $\gamma$
factors through
$$\Tor_{\bZ}\bigl(\Om_*(B\pi'),\Om_*(B\bZ/p^k)\bigr).$$ But this
latter group is generated by cobordism Massey products $\langle M',
p^l, M'' \rangle$, where $M'\to B\pi'$ and $M''\to B\bZ/p^k$. By
inductive hypothesis, $M'$ and $M''$ can be chosen to have {\nYi},
hence so can this Toda bracket by Theorem \ref{thm:nnYiToda}. Again,
the exceptional case of that Theorem never arises in our context. So
this shows that the image of $\gamma$ is represented by manifolds of
{\nYi}. Chasing diagram \eqref{eq:Torchase} now shows that
$RH_*(B\pi)$ is represented by manifolds with {\nYi}, which completes
the inductive step for the case $p$ odd.
\end{labeleddemo}
\begin{cor}\label{cor:abelianSylow}
Let $M^n$ be a closed, connected, $n$-manifold with finite
fundamental group $\pi$, with non-spin universal cover, and with $n\ge
5$. Also assume all Sylow subgroups of $\pi$ are abelian.
Then $M$ has \nYi.
\end{cor}
\begin{demo}This is immediate from Theorem \ref{thm:nnYinonspinabel}
and Theorem \ref{thm:KS}.
\end{demo}
In the odd order case, we can carry this over to the spin case as well:
\begin{thm}\label{thm:oddspinabelianSylow}
Let $M^n$ be a closed, connected, spin $n$-manifold with finite
fundamental group $\pi$ of odd order, and with $n\ge
5$. Also assume all Sylow subgroups of $\pi$ are abelian.
Then $M$ has \nYi.
\end{thm}
\begin{demo}By Petean's theorem \cite{P2}, this is true when $\pi$ is
trivial. As before, it's enough to consider the case of an abelian
$p$-group, $p$ odd. But for $\pi$ of odd order, the natural map $\widetilde
\Oms_*(B\pi) \to \widetilde\Om_*(B\pi)$ is an isomorphism, since the
map of spectra $M\text{Spin}\to MSO$ is an equivalence after
localizing at $p$ (see \cite{MM}). We prove the result by induction on
the rank of $\pi$. Thus write $\pi=\pi'\times \bZ/p^k$, and assume by
inductive hypothesis that the theorem is true for $\pi'$. Since 
$\widetilde \Oms_*(B\pi') \to \widetilde\Om_*(B\pi')$ and
$\widetilde \Oms_*(B\bZ/p^k) \to \widetilde\Om_*(B\bZ/p^k)$ are
isomorphisms, we have by Landweber's Theorem \cite{Landw} an exact
sequence
\begin{multline}\label{eq:spinTorseq}
0\to   \Oms_*(B\pi')\otimes_{\Oms_*} \Oms_*(B\bZ/p^k)
\to  \Oms_*(B\pi) \\ \to  
\Tor_{\Oms_*}\bigl(\Oms_*(B\pi'),\Oms_*(B\bZ/p^k)\bigr) \to  0.
\end{multline}
By inductive hypothesis, each element of $\Oms_s(B\pi')$ is
represented by a map $M'\to B\pi'$, with $M'$ a spin $s$-manifold with
{\nYi}, and similarly  each element of $\Oms_t(B\bZ/p^k)$ is
represented by a map $M''\to B\bZ/p^k$, with $M''$ a spin $t$-manifold with
{\nYi}. Then $[M'\to B\pi']\otimes [M''\to B\bZ/p^k]$ in
the tensor term on the left side of \eqref{eq:spinTorseq}
is represented by $M'\times M''\to B\pi$, which has {\nYi}.
Furthermore, the Tor term
$\Tor_{\Oms_*}([M'\to B\pi'], [M''\to B\bZ/p^k])$ on the right
of \eqref{eq:spinTorseq} pulls back (non-canonically) to the class of
a Toda bracket $\langle M', P, M''\rangle$, where $P$ is some element
from the coefficient ring $\Oms_*$. Since $M'$, $P$, and $M''$ all have
{\nYi} (we don't even need to know anything about $P$!), it follows
from Theorem \ref{thm:nnYiToda} that this Toda bracket has {\nYi}, and
this completes the inductive step.
\end{demo}

\section{Applications to Positive Scalar Curvature}\label{sec:psc}

It turns out that the method of proof of Theorem \ref{thm:nnYinonspinabel},
if we replace Theorem \ref{thm:nnYiToda} by Theorem \ref{thm:pscToda}, gives
substantial results on the {\psc} problem for manifolds
with finite abelian fundamental group for which the universal
cover is non-spin, since all of the homology generators
constructed above have {\psc} by
Theorem \ref{thm:pscToda}, \emph{except} for those involving
Toda brackets and products of one-dimensional homology classes.
We proceed to formalize this as follows:
\begin{defn}\label{def:toral}
Let $\pi$ be a finitely generated abelian group.
Call a class in $H_n(B\pi,\bZ)$ or in $H_n(B\pi,\bZ/2)$
\emph{toral} if it is represented
by a map $T^n\to B\pi$. Note that any such map is determined
up to homotopy by the associated map $\bZ^n\to \pi$ on fundamental 
groups, which we may assume without loss of generality
to have image of rank $n\le r$,
where $r$ is the \emph{rank} of $\pi$, that is, the minimal number of
cyclic factors when we write $\pi$ as a direct sum of cyclic
groups, so toral classes only exist in degrees $n\le r$. 
Let $H_n^{\text{toral}}(B\pi,\bZ)\subseteq H_n(B\pi,\bZ)$ be the 
subgroup generated by the the toral classes, and call this
the \emph{toral subgroup}.
\end{defn}
\begin{prop}\label{prop:toral}
Let $\pi$ be an elementary abelian $p$-group of rank $r$, that is, 
$(\bZ/p)^r$. Then for all $n\ge 1$, $H_n(B\pi,\bZ)$ is also
elementary abelian, of rank equal to 
\[ \sum_{j=1}^n(-1)^{n-j}\left(
\begin{matrix}j+r-1\\r-1\end{matrix}\right).
\]
The toral subgroup $H_n^{\mathrm{toral}}(B\pi,\bZ)$ is a direct
summand in $H_n(B\pi,\bZ)$, of rank the binomial coefficient
$\left(\begin{matrix}r\\n\end{matrix}\right)$. {\lp}Note that this
vanishes for $n>r$.{\rp}
\end{prop}
\begin{demo}The homology groups $H_n(B\bZ/p,\bZ)$ vanish for
$n>0$ even and are $\bZ/p$ for $n$ odd. So
by iterated applications of the K\"unneth Theorem,
all integral homology groups of $\pi$ (other than $H_0$, which is
$\bZ$), are elementary abelian $p$-groups.
Consider the Poincar\'e series 
\[
P(r,t)=1+\sum_{n=1}^\infty t^n \dim_{\bZ/p} H_n(B(\bZ/p)^r,\bZ).
\]
Then 
\begin{equation}\label{eq:P1t}
P(1,t)=1+t+t^3+t^5+\cdots=1+\frac{t}{1-t^2}=\frac{1+t-t^2}{1-t^2}.
\end{equation}
The K\"unneth Theorem gives the recursion relation
\begin{equation}\label{eq:Prt}
P(r+1,t)=P(r,t)P(1,t)+t\left(P(r,t)-1\right)\left(P(1,t)-1\right),
\end{equation}
where the first term comes from the ``tensor terms'' and the
second term comes from the ``Tor terms.'' Putting together
equations (\ref{eq:P1t}) and (\ref{eq:Prt})
yields by induction on $r$ the formula
\[
P(r,t)=\frac{1+t(1-t)^r}{(1-t)^r(1+t)}
=\frac{t}{1+t}+\frac{1}{(1-t)^r(1+t)}.
\]
For $n\ge 1$, the coefficient of $t^n$ in this expression is
\[
(-1)^{n+1}+\sum_{j=0}^n(-1)^{n-j}\left(
\begin{matrix}j+r-1\\r-1\end{matrix}\right),
\]
which is the expression in the statement of the Proposition.
On the other hand, the toral subgroup is generated just by the
products of distinct generators of $H_1$, so in degree $n$,
we have $\left(\begin{matrix}r\\n\end{matrix}\right)$
possibilities.
\end{demo}
Similarly (though more simply), we have:
\begin{prop}\label{prop:toralnonorient}
Let $\pi$ be an elementary abelian $2$-group of rank $r$, that is, 
$(\bZ/2)^r$. Then for all $n\ge 0$, $H_n(B\pi,\bZ/2)$ is also
elementary abelian, of rank equal to 
\[ \left(\begin{matrix}n+r-1\\r-1\end{matrix}\right).
\]
The toral subgroup $H_n^{\mathrm{toral}}(B\pi,\bZ/2)$ is a direct
summand in $H_n(B\pi,\bZ/2)$, of rank the binomial coefficient
$\left(\begin{matrix}r\\n\end{matrix}\right)$. {\lp}Note that this
vanishes for $n>r$.{\rp}
\end{prop}
\begin{demo}This is easier than the previous case since there
are no Tor terms.  The Poincar\'e series for $H_*(B\bZ/2,\bZ/2)$
is $\sum_{j=0}^\infty t^n = 1/(1-t)$, and so the
Poincar\'e series for $H_*(B\pi,\bZ/2)$ is 
\[
\frac{1}{(1-t)^r}=\sum_{n=0}^\infty 
\left(\begin{matrix}n+r-1\\r-1\end{matrix}\right) t^n,
\]
which yields the desired formula.
\end{demo}
\begin{defn}\label{def:atoral}
Let $\pi$ be a finite abelian $p$-group, say $\prod_{i=1}^r
\bZ/p^{k_i}$. For purposes of this definition, $H_*$ means homology
with $\bZ$ coefficients if $p\ne 2$, and homology
with $\bZ/2$ coefficients if $p= 2$. We define a splitting
$H_n(B\pi)= H_n\tor(B\pi) \oplus H_n\ator(B\pi)$ inductively as
follows. If $r=1$, $H_n\tor(B\pi)=H_n(B\pi)$ if $n\le 1$
and $H_n\tor(B\pi)=0$ for $n\ge 2$, so we let
$H_*\ator(B\pi)=\bigoplus_{n\ge 2}H_n(B\pi)$. If $r>1$,
write $\pi=\pi'\times \bZ/p^k$, $k=k_r$, where
$\pi'=\prod_{i=1}^{r-1} \bZ/p^{k_i}$. Choose a splitting of the
K\"unneth formula
\begin{equation}\label{eq:Kun}
0\to   H_*(B\pi')\otimes_{\bZ} H_*(B\bZ/p^k)
\to  H_*(B\pi) \leftrightarrows  
\Tor_{\bZ}\bigl(H_*(B\pi'),H_*(B\bZ/p^k)\bigr) \to  0.
\end{equation}
By induction, we have splittings $H_*(B\pi')=H_*\tor(B\pi')\oplus
H_*\ator(B\pi')$ and $H_*(B\bZ/p^k)=H_*\tor(B\bZ/p^k)\oplus
H_*\ator(B\bZ/p^k)$. The atoral part of the tensor term on the left is
defined to be 
\[
\bigl(H_*(B\pi')\otimes H_*\ator(B\bZ/p^k)\bigr) \oplus
\bigl(H_*\ator(B\pi')\otimes H_*\tor(B\bZ/p^k)\bigr).
\]
Then we define $H_*\ator(B\pi)$ to be the direct sum of the atoral piece of
the tensor term with the image of the $\Tor$ term under
the splitting of the exact sequence \eqref{eq:Kun}.
\end{defn}
\begin{defn}\label{def:RH}
For any space $X$, we denote by $RH_*(X)$ the image of
the Hurewicz map
$\Om_*(X)$ in $H_*(X,\bZ)$, and call it the \emph{representable homology}.
(This already made an appearance in the proof of Theorem
\ref{thm:nnYinonspinabel}.) Note that $RH_*$ is a functor, though not
a homology theory. By Lemma \ref{lem:cyclic}, $RH_*(B\pi)=H_*(B\pi,\bZ)$
when $\pi$ is a cyclic group.
\end{defn}
The following fact, which is somewhat surprising, will be our key
technical tool:
\begin{prop}\label{prop:ranktwo}
Let $\pi$ be an elementary abelian $p$-group of rank $2$, where $p$ is an odd
prime. Then $RH_{\text{\textup{odd}}}(B\pi)$ is generated 
{\lp}as an abelian group{\rp} by
the images of $RH_*(B\sigma)$, as $\sigma$ runs over the cyclic
subgroups of $\pi$. 
\end{prop}
\begin{demo}This is
proved in \cite{BG}, using explicit calculations of the eta-invariants
of lens spaces. \end{demo}
\begin{thm}\label{thm:RHstructure}
Let $\pi$ be an elementary abelian $p$-group, where $p$ is an odd
prime. Then $RH_*(B\pi)$ is generated {\lp}as an abelian group{\rp}
by elements $$x_1\otimes \cdots \otimes x_j \in
H_*(B \sigma_1)\otimes \cdots \otimes H_*(B \sigma_j),$$
with $\sigma_1\times \cdots \times \sigma_j$ a subgroup of $\pi$
with each $\sigma_i$ a cyclic $p$-group.
\end{thm}
\begin{demo}We prove this by induction on the rank $r$. When $r=1$,
the statement is trivially true, and when $r=2$, this is
Proposition \ref{prop:ranktwo}. Now assume the result for smaller
values of $r$, and write $\pi=\pi'\times\bZ/p$, where $\pi'$ has rank
$r-1$. We again use diagram \eqref{eq:Torchase} (with $k=1$).
The image of $\alpha$ is taken care of by inductive hypothesis. So
consider the image of $\gamma$. Consider a representable
class $\Tor(x,y)$, where $y\in
H_*(B\bZ/p)$ and $x\in RH_*(B\pi')$. By inductive hypothesis, we may
assume $x$ is the image of a representable
class $x_1\otimes \cdots \otimes x_j \in
H_*(B \sigma_1)\otimes \cdots \otimes H_*(B \sigma_j)$,
with $\sigma_1\times \cdots \times \sigma_j$ a subgroup of $\pi'$
with each $\sigma_i$ a cyclic $p$-group. It will suffice to show that
$\Tor(x_1\otimes \cdots \otimes x_j , y)$ is of the correct form back
in $H_*(B(\sigma_1\times \cdots \times \sigma_j\times
\bZ/p))$. Represent each $x_i$ by a manifold $M_i\to B\sigma_i$ which
is either $S^1$ or a lens space, and also represent $y$ by a manifold
$L$ which is either $S^1$ or a lens space. Then $\Tor(x_1\otimes
\cdots \otimes x_j , y)$ is represented by the homology Toda bracket $\langle
[M_1]\times \cdots \times [M_j], p, [L]\rangle$. Since this Toda
bracket is representable, it must be that $[M_1]\times \cdots \times
[M_j]$ and $[L]$ have order $p$ in bordism, so that the homology Toda
bracket lifts to a bordism Toda bracket, giving a class in
$\Om_*(B(\sigma_1\times \cdots \times \sigma_j\times
\bZ/p))$. Indeed, the representable part of the Tor term in homology
is the group generated by the
images of the bordism Toda brackets $\langle [M], p^s, [L]\rangle$, 
and such a bracket maps to zero in homology whenever $s>1$.
Now consider the bordism Toda bracket $\langle
[M_1]\times \cdots \times [M_j], p, [L]\rangle$.
We know at least one of the $[M_i]$, say $[M_j]$, has order $p$.
Then
\[
\langle [M_1]\times \cdots \times [M_j], p, [L]\rangle
= [M_1]\times \cdots \times [M_{j-1}]\times \langle [M_j], p,
[L]\rangle
\]
(\cite{Al}, 2.1). Now apply Proposition \ref{prop:ranktwo}
to $\langle [M_j], p, [L]\rangle$, and this completes the inductive step.
\end{demo}
\begin{thm}\label{thm:pscfiniteabelian}
Let $\pi$ be an abelian $p$-group, and let $n\ge 5$.
Suppose $M^n$ is a closed manifold with fundamental group $\pi$ and
non-spin universal cover, and suppose $[M\to B\pi]\in
H_n\ator(B\pi)$. If $p=2$, also assume $M$ is not orientable.
If $p$ is odd, also assume $\pi$ is elementary abelian.
Then $M$ has a metric with \psc.
In particular, if $n> \operatorname{rank}\pi$, then
\textbf{every} $n$-manifold with fundamental group $\pi$ and
with non-spin universal cover has a metric of \psc.
\end{thm}
\begin{demo}First consider the non-orientable case with $p=2$. In this
case, since there are no Tor terms,
\[ 
H_n(B\pi,\bZ/2) =\bigotimes_{i=1}^{\operatorname{rank}\pi} 
H_*(B\pi_i,\bZ/2),
\]
where the $\pi_i$ are the cyclic factors of $\pi$. 
Now each class in $H_*(B\pi_i,\bZ/2)$ may by Lemma
\ref{lem:cyclic} be represented by a
manifold with \nYi, or with {\psc} if the class is in degree $>1$. 
So it immediately follows that $H_n\ator(B\pi,\bZ/2)$ is
represented by manifolds of \psc. That takes care of the
non-orientable case.

Now consider the case where $p$ is odd (and $M$ is orientable).
We apply Theorem \ref{thm:RHstructure}. This reduces us to the
case of classes of the form $x_1\otimes \cdots \otimes x_j \in
H_*(B \sigma_1)\otimes \cdots \otimes H_*(B \sigma_j)$,
with each $x_i$ represented by either $S^1$ or a lens space.
The product manifold has {\psc} unless all the $x_i$'s are
$1$-dimensional, in which case we have a toral class.
\end{demo}

\begin{prob}Are toral homology classes {\lp}for an elementary abelian
$p$-group{\rp} represented by manifolds of \psc? We suspect not,
but we know of no way to approach this question.
\end{prob}

\begin{prob}Is Theorem \ref{thm:oddspinabelianSylow} true without the
odd order assumption?  We presume so, but the proof would necessarily
be much more complicated, since computing $ko_*(B\pi)$ for a $2$-group
is quite difficult.
\end{prob}

\begin{prob}Is Theorem \ref{thm:pscfiniteabelian} true for arbitrary
abelian $p$-groups? Again we suspect so, but the necessary
calculations are difficult.
\end{prob}

\begin{flushleft}
\textbf{Author Addresses:}\\
\vspace*{.2in}
Boris Botvinnik\\
Department of Mathematics\\
University of Oregon\\
Eugene, OR 97403-1222\\
\texttt{botvinn@poincare.uoregon.edu}\\
\vspace*{.2in}
Jonathan Rosenberg\\
Department of Mathematics\\
University of Maryland\\
College Park, MD 20742-4015\\
\texttt{jmr@math.umd.edu}
\end{flushleft}

\begin{thebibliography}{99}

\bibitem{Al}J.\ C.\ Alexander, 
Cobordism Massey products, 
Trans.\ Amer.\ Math.\ Soc.\ \textbf{166} (1972), 197--214. 
\bibitem{And}M.\ T.\ Anderson, 
Scalar curvature, metric degenerations and the static vacuum 
Einstein equations on $3$-manifolds, I, 
Geom.\ Funct.\ Anal.\ \textbf{9} (1999), no.\ 5, 855--967. 
\bibitem{BG}B.\ Botvinnik and  P.\ Gilkey, 
The eta invariant and the Gromov-Lawson conjecture for
elementary abelian groups of odd order, Topology Appl.\ \textbf{80} (1997),
no.\ 1--2, 43--53. 
\bibitem{BGS}B.\ Botvinnik, P.\ Gilkey, and S.\ Stolz,
The Gromov-Lawson-Rosenberg conjecture for groups with periodic cohomology, 
J.\ Differential Geom.\ \textbf{46} (1997), no.\ 3, 374--405.
\bibitem{GL1}M.\ Gromov and H.\ B.\ Lawson, Jr.,
Spin and scalar curvature in the presence of a fundamental group, I,
Ann.\ of Math.\ \textbf{111}  (1980),   209--230.   
\bibitem{GL2}M.\ Gromov and H.\ B.\ Lawson, Jr.,
The classification of simply connected manifolds of positive scalar
curvature, Ann.\ of Math.\
\textbf{111}  (1980),  423--434.                     
\bibitem{GL3}M.\ Gromov and H.\ B.\ Lawson, Jr.,
Positive scalar curvature and the Dirac operator on complete
Riemannian manifolds, Publ.\ Math.\ I.H.E.S., no.\  58
(1983),  83--196.                                      
\bibitem{JW}
D.\ C.\ Johnson and S.\ W.\ Wilson, 
The Brown-Peterson homology of elementary $p$-groups,
Amer.\ J.\ Math.\ \textbf{107} (1985), no.\ 2, 427--453. 
\bibitem{KS1}S.\ Kawsik and R.\ Schultz, Positive scalar curvature 
and periodic fundamental
groups, Comment.\ Math.\ Helv.\ \textbf{65} (1990), no.\ 2, 271--286. 
\bibitem{KW} J.\ L.\ Kazdan and F.\ Warner,
Existence and
   conformal deformation of metrics with prescribed Gaussian and scalar
   curvatures, Ann.\ of Math.\ (2) \textbf{101} (1975), 317--331.
\bibitem{Kob}O.\ Kobayashi, 
Scalar curvature of a metric with unit volume,
Math.\ Ann.\ \textbf{279} (1987), no.\ 2, 253--265.
\bibitem{Landw}P.\ Landweber, 
K\"unneth formulas for bordism theories,
Trans.\ Amer.\ Math.\ Soc.\ \textbf{121} (1966), no.\ 1, 242--256. 
\bibitem{LeB}C.\ LeBrun, Kodaira dimension and the Yamabe problem, 
Comm.\ Anal.\ Geom.\ \textbf{7} (1999),
no.\ 1, 133--156. 
\bibitem{Loh}J.\ Lohkamp, The space of negative {\scurv} metrics,
Invent.\ Math.\ \textbf{110} (1992), 403--407.
\bibitem{MM}I.\  Madsen and R.\ J.\ Milgram, 
\emph{The classifying spaces for surgery and cobordism of manifolds},
Ann.\ of Math.\ Studies, vol.\ 92,
Princeton Univ.\ Press, Princeton, N.J., 1979. 
\bibitem{P1}J.\ Petean, Computations of the Yamabe invariant, 
Math.\ Res.\ Lett.\ \textbf{5} (1998), no.\ 6,
703--709. 
\bibitem{P2}J.\ Petean,  The Yamabe invariant of simply connected
manifolds, J.\ reine und angew.\ Math.\ \textbf{523} (2000), 225--231.
\bibitem{PY}J.\ Petean and G.\ Yun, Surgery and the Yamabe invariant,
Geom.\ and Funct.\ Anal.\ \textbf{9} (1999), 1189--1199.
\bibitem{RosKO}J.\ Rosenberg,
The $KO$-assembly map and positive scalar
curvature, in \emph{Algebraic Topology Pozna\'n 1989}, Lecture Notes
in Math., vol.\ 1474, Springer, Berlin, 1991, pp.\ 170--182.
\bibitem{RosCob}J.\ Rosenberg,
Reflections on C.\ T.\ C.\ Wall's work on cobordism, 
in \emph{Surveys on Surgery Theory: Volume 2}, S.\ Cappell,
A.\ Ranicki, and J.\ Rosenberg, eds., Annals of Math.\ Studies,
vol.\ 149, Princeton Univ.\ Press, Princeton, NJ, 2001, pp.\ 49--61.
\bibitem{RS1}J.\ Rosenberg and S.\ Stolz,
A ``stable'' version of the Gromov-Lawson conjecture, in
\textit{The \v Cech Centennial: Proc.\ Conf. on Homotopy Theory},
M.\ Cenkl and H.\ Miller, eds., Contemp.\ Math.\ \textbf{181} (1995),
Amer.\ Math.\ Soc.,  Providence, RI, pp.\ 405--418.
\bibitem{RS2}J.\ Rosenberg and S.\ Stolz,
Metrics of \psc\ and connections with surgery, 
in \emph{Surveys on Surgery Theory: Volume 2}, S.\ Cappell,
A.\ Ranicki, and J.\ Rosenberg, eds., Annals of Math.\ Studies,
vol.\ 149, Princeton Univ.\ Press, Princeton, NJ, 2001, pp.\ 353--386.
\bibitem{Sc}R.\ Schoen,
Conformal deformation of a Riemannian metric to constant scalar curvature,
J.\ Differential Geom.\ \textbf{20} (1984), no.\ 2, 479--495.
\bibitem{SY}R.\ Schoen and S.-T.\ Yau,
On the structure of manifolds with positive scalar curvature.
Manuscripta Math.\ \textbf{28} (1979), no.\ 1-3, 159--183.
 \bibitem{Sch}R.\ Schultz, Positive scalar curvature and 
odd order abelian fundamental groups,
Proc.\ Amer.\ Math.\ Soc.\ \textbf{125} (1997), no.\ 3, 907--915. 
\bibitem{S1}S.\ Stolz,
Simply connected manifolds of positive scalar curvature,
Ann.\ of Math.\ (2) \textbf{136} (1992), no.\ 3, 511--540. 
\bibitem{S2}S.\ Stolz, Positive scalar curvature 
metrics---existence and classification questions, in
\emph{Proceedings of the International Congress of Mathematicians},
Vol.\ 1 (Z\"urich, 1994), Birkh\"auser,
Basel, 1995, pp.\ 625--636. 

\bibitem{Stong} R.\ Stong, \textit{Notes on Cobordism Theory},
Mathematical Notes, Princeton Univ.\ Press, Princeton, NJ, 1968.
\bibitem{Tay}L.\ R.\ Taylor, 
$2$-local cobordism theories,
J.\ London Math.\ Soc.\ (2) \textbf{14} (1976), no.\ 2, 303--308. 
\end{thebibliography}
\end{document}